\documentclass[12pt]{amsart}

\usepackage{amsfonts,amsmath,amsthm}
\usepackage{latexsym}

\newtheorem{theorem}{Theorem}[section]
\newtheorem{corollary}[theorem]{Corollary}
\newtheorem{lemma}[theorem]{Lemma}

\newtheorem{definition}{Definition}[section]
\theoremstyle{definition}

\newcommand \tr {\mathrm{Tr}\:}
\newcommand \im {\mathrm{Im}\:}

\newcommand{\C}{\ensuremath{\mathbb{C}}}
\newcommand{\R}{\ensuremath{\mathbb{R}}}
\newcommand{\K}{\ensuremath{\mathbb{K}}}
\newcommand{\adj}{^\ast}
\newcommand{\sq}{^{1/2}}

\def \e {\varepsilon}

\def \l {\lambda}

\def \< {\left\langle}
\def \> {\right\rangle}

\begin{document}

\title[On the linear polarization constant of $\R^n$]{Linear polarization constant of $\R^n$}
\author{M\'at\'e~Matolcsi}
\address{ Alfr\'ed R\'enyi Institute of Mathematics,
Hungarian Academy of Sciences POB 127 H-1364 Budapest, Hungary
Tel: (+361) 483-8302, Fax: (+361) 483-8333}
\email{matomate@renyi.hu}
%%% \date{\today}

\begin{abstract}
The present work contributes to the determination of the $n$-th
linear polarization constant $c_n(H)$ of an $n$-dimensional real
Hilbert space $H$. We provide some new lower bounds on the value
of $\sup_{\|y\|=1}|\< x_1,y \> \cdots \< x_n,y \> |$, where $x_1,
\dots ,x_n$ are unit vectors in $H$. In particular, the results
improve an earlier estimate of Marcus. However, the intriguing
conjecture $c_n(H)=n^{n/2}$ remains open.
\end{abstract}

\maketitle



\section{Introduction}
In this note we aim to make a contribution to estimating the $n$-th linear
polarization constant $c_n(H)$ of an $n$-dimensional real Hilbert space $H$. We
begin with introducing some (more general) standard terminology
and giving a short account of some related results.

Let $X$ denote a Banach space over the real or complex field $\K$. A function
$P: \ X\to \K$ is a continuous {\it n-homogeneous polynomial} if there exists
a continuous $n$-linear form $L: \ X^n\to\K$ such that
$P(x)=L(x,\dots ,x)$ for all $x\in X$. We define
$$\|P\|:=\sup \{ |P(x)|: \ x\in B\}$$
where $B$ denotes the unit ball of $X$. Considerable attention has been
devoted to polynomials of the form $P(x)=f_1(x)f_2(x)\dots f_n(x)$, where
$f_1, f_2,\dots ,f_n$ are bounded linear functionals on $X$. For any{\it complex}
Banach space $X$ Ben\'\i tez, Sarantopoulos and Tonge \cite{ben} have obtained
$$\|f_1\| \ \|f_2\|\cdots \|f_n\|\le n^n\|f_1 f_2 \dots f_n\|,$$
and the constant $n^n$ is best possible. For {\it real} Banach spaces, Ball's
solution \cite{ball2} of the famous plank problem of Tarski gives the same
general result. For specific spaces, however, the general constant $n^n$ can
be lowered.

\begin{definition} (Ben\'\i tez, Sarantopoulos, Tonge \cite{ben})
The n-th linear polarization constant of a Banach space $X$ is defined by
\begin{eqnarray*}\label{polconstdef}
c_n(X)&:=& \inf \{M : \|f_{1}\| \cdots \|f_{n}\| \leq M
\|f_{1}\cdots f_{n}\|\; (\forall f_{1}, \ldots ,f_{n} \in
X^{\ast}) \}
 \\
 & = & 1/ \inf_{f_{1}, \ldots ,f_{n} \in S_{X^{\ast}}}
\sup_{\|x\|=1}|f_{1}(x) \cdots f_{n}(x)|.
\end{eqnarray*}
The linear polarization constant of $X$ is defined by
\begin{equation}\label{limitconstant}
c(X):=\lim_{n \rightarrow \infty} c_{n} (X)^{ \frac{1}{n}}\;\;.
\end{equation}
\end{definition}

Let us recall that the above definition of $c(X)$ is justified since
R\'ev\'esz and Sarantopoulos \cite{sar} showed that the limit
\eqref{limitconstant} does exist. Moreover, they also showed (both
in the real and complex cases) that $c(X)=\infty$ if and only if
$\dim X=\infty$.

Note that it is easy to see that for any Banach space $X$ we have
\begin{equation}\label{finite}
c_{n}(X)=\sup \left\{c_{n} (Y): Y {\mbox {is a closed subspace
of}}\, X, \, \dim Y=n \right\}\,.
\end{equation}
In particular, for a real or complex Hilbert space $H$
of dimension at least $n$, we
always have $c_n(H)=c_n(\K^n)$.

Ben\'\i tez, Sarantopoulos and Tonge \cite{ben} proved that for
isomorphic Banach spaces $X$ and $Y$ we have $c_n(X)\le
d^n(X,Y)c_n(Y)$, where $d(X,Y)$ denotes the Banach-Mazur distance
of $X$ and $Y$. Note, that for any $n$-dimensional space $X$ a
result of John \cite{john} states that $d(X,\K^n)\le \sqrt{n}$
(where $\K^n$ denotes the $n$-dimensional Hilbert space).  The
combination of these results mean that the determination of
$c_n(\K^n)$ gives information on the linear polarization constants
of other spaces, too.

Here we are going to focus our attention to Hilbert spaces.
Pappas and R\'ev\'esz \cite{papp} showed that $c(\K^n)= e^{-L(n,\K
)}$, where
$$L(n,\K ):= \int_{S}\textrm{log}|\< x,e\>  |d\sigma (x);$$
here $S$ and $\sigma$  denote the unit sphere and the normalized
surface measure, respectively, and $e\in S$ is an arbitrary unit
vector. This result gives information on the asymptotic behaviour
of $c_m(\K^n)$ as $m\to \infty$.
However, the exact values of $c_m(\K^n)$ seem very hard to determine.
Anagnostopoulos and R\'ev\'esz \cite{anag} found explicit
relations between $c_n(\R^2)$, $c_n(\C^2)$ and the Chebyshev
constants of $S^1$ and $S^2$, respectively. These relations, however, 
do not seem to carry over to higher dimensional spaces. Note that
the $n$-th Chebyshev constant of $S^1$ is well-known, but the
exact determination of the $n$-th Chebyshev constant of $S^2$ seems
hopeless; see \cite{anag} and \cite{wag}. 

A remarkable result of
Arias-de-Reyna \cite{rey} states that $c_n(\C^n)=n^{n/2}$. Ball's
recent solution \cite{ball} of the complex plank problem also
implies the same result.

The value of $c_n(\R^n)$ seems even harder to find. The
determination of $c_n(\R^n)$, by the definition and the Riesz
representation theorem, boils down to determining
$$
I:=
\inf_{x_{1}, \dots ,x_{n} \in S} \sup_{\|y\|=1}|\< x_{1},y\>
\cdots \< x_{n},y\>  |
$$
The estimate $I\le n^{-\frac{n}{2}}$ follows by considering an
orthonormal system.

The result of Arias-de-Reyna can be used to derive the following
estimates \cite{sar}:
$$
n^{\frac{n}{2}}\le c_n(\R^n)\le 2^{\frac{n}{2}-1}n^{\frac{n}{2}}.
$$
A natural, intriguing conjecture, see \cite{ben}, \cite{sar} is
the following.

\medskip
\noindent {\bf Conjecture.} $c_n(\R^n)=n^{n/2}$.
\medskip

Marcus (communicated in \cite{marc}, and elaborated later in
\cite{sar}) gives the following estimate: If $x_1, x_2, \dots ,
x_n$ are unit vectors in $\R^n$ then there exists a unit vector
$y$ such that
\begin{equation}\label{marcus}
|\< x_1,y\> \cdots\< x_n,y\> |\ge (\l_1/n)^{n/2},
\end{equation}
where $\l_1$ denotes the smallest eigenvalue of the Gram matrix
$XX\adj =[\< x_i,x_j\>  ]$. Marcus also expressed the opinion
that lower bounds on $\sup_{\|y\|=1}|\< x_1,y\>  \cdots \<
x_n,y\>  |$ should involve
%%% a 'measure of orthogonality' of the vectors $x_1, x_2, \dots
%%% ,x_n$;
%%% $\sup_{\|y\|=1}|\< x_1,y\>  \cdots \< x_n,y\>  |\ge f(\l_1, \dots
%%% , \l_n) n^{-\frac{n}{2}}$, where $\l_1, \dots , \l_n$ are
the eigenvalues $\l_1,\dots,\l_n$ of the Gram matrix  $XX\adj
=[\< x_i,x_j\>  ]$, i.e. we should look for estimates of the form
$\sup_{\|y\|=1}|\< x_1,y\>  \cdots \< x_n,y\>  |\ge f(\l_1, \dots
, \l_n) n^{-n/2}$. Note that $\sum_j \l_j=\tr XX\adj=n$.
Therefore the above Conjecture can be equivalently formulated as
$$
\sup_{\|y\|=1}|\< x_1,y\>  \cdots \< x_n,y\>  | \ge 1\cdot
n^{-n/2}= \left(\frac{\l_1 +\dots +\l_n}{n}\right)^{n/2}n^{-n/2}.
$$

In the next section we show that if we replace the arithmetic mean
by the harmonic mean of the numbers $\l_1, \dots , \l_n$, then
the corresponding lower estimate does hold. This gives an
improvement of Marcus' result but the Conjecture still remains
open.

\section{Lower bounds}

In this section we are going to present three different lower
bounds on the value of $\sup_{\|y\|=1}|\< x_1,y\>  \cdots \<
x_n,y\>  |$. The first result relies on an averaging argument,
while the next two uses the following lemma of Bang \cite{bang}:

\begin{lemma}
Let $H=(h_{jk})$ be an $n\times n$ Gram matrix and $r_1, \dots r_n$ be a sequence of positive numbers. Then there are signs $\e_1, \dots \e_n$ for which
\begin{equation}\label{bang}
\e_jr_j\sum_{k=1}^n h_{jk}r_k\e_k \ge r_j^2
\end{equation}
for every $j$.
\end{lemma}

Assume now, that we are given $n$ unit vectors $x_1, \dots ,x_n$
in $\R^n$. We are looking for a unit vector $y$ such that the
product $|\< x_1,y\>  \cdots \< x_n,y\>  |$ is 'as large as
possible'. Let $X$ denote the $n\times n$ matrix whose $j$-th row
is $x_j$. With this notation, our aim is to maximize the
expression $\prod_{j=1}^{n}|(Xy)_j|$. As a first observation we
'symmetrize' the matrix $X$. Take the polar decomposition
$X\adj =V(XX\adj )\sq$ of $X\adj$. The partial
isometry $V$ maps $(\ker X\adj )^\bot$ to $\im X\adj$, and can be
extended (not uniquely, in general) to a unitary operator $U$.
Taking adjoints we get $X=(XX\adj )\sq U\adj$. The unitary
operator $U\adj$ maps the unit sphere onto itself, therefore
$$
\max_{\|y\|=1}\prod_{j=1}^{n}|(Xy)_j|= \max_{\|y\|=1}\prod_{j=1}^{n}|((XX\adj )\sq y)_j|.
$$
Therefore, in the forthcoming arguments, we are going to consider the positive
self-adjoint matrix $(XX\adj )\sq$ instead of the original matrix $X$. Note, also, that the polar decomposititon implies that the rows (and , by symmetry, the columns) of the matrix  $(XX\adj )\sq$ are also unit vectors

\begin{theorem}\label{prop1}
Assume that the given unit vectors $x_1,\dots , x_n$ are linearly independent.
Let $V_1, \dots ,V_n$ denote the lengths of the column-vectors of
the matrix $(XX\adj )^{-1/2}$. Then
\begin{equation}\label{opt}
\max_{\|y\|=1}\prod_{j=1}^{n}|(Xy)_j|
\ge \frac{1}{V_1 \cdots V_n}\cdot n^{-\frac{n}{2}}
\end{equation}
\end{theorem}

\begin{proof}
Note that the inverse matrix $(XX\adj )^{-1/2}$ exists by the
assumption of linear independency. Denote the entries of
$(XX\adj )^{-1/2}$ by $(v_{jk})$.

Take any vector $c:=(c_1,\dots ,c_n)^T$ with $c_j\ge0 \,\,
(j=1,\dots,n)$ and consider all possible $2^n$ arrangements of
signs of the entries in $c$: i.e., consider $ c(\e):=(\e_1c_1,
\dots , \e_nc_n)^T$ for all $\e\in \{\pm 1\}^n$. We evaluate the
sum $\sum_\e \|(XX\adj )^{-1/2} c (\e ) \|^2$. It is easy to see
that the double products $2c_jv_{ij}c_kv_{ik}$ all cancel out
because they appear an equal number of times with positive and
negative signs. Therefore,
\begin{equation}
\sum_\e \|(XX\adj )^{-1/2} c(\e ) \|^2 = 2^n \cdot
\sum_{jk}c_k^2v_{jk}^2 =2^n \cdot \sum_k c_k^2V_k^2
\end{equation}

Hence, there exists a signed vector  $c(\e ) =(\e_1c_1, \dots ,
\e_nc_n)^T$ such that $\|(XX\adj )^{-1/2} c(\e) \|\le (\sum_k
c_k^2V_k^2)\sq $. Take
$$y:=\frac{(XX\adj )^{-1/2} c(\e) }{\|(XX\adj )^{-1/2} c(\e)\|}.$$
Then
\begin{equation}\label{becs}
\prod_{j=1}^{n}|((XX\adj )\sq y)_j|
\ge \prod_j c_j \cdot (\sum_k c_k^2V_k^2)^{-n/2}.
\end{equation}

After introducing the new variables $b_1:=c_1V_1, \dots
b_n:=c_nV_n$, the right hand side becomes $\prod_j b_j \cdot
(\sum_k b_k^2)^{-n/2} \cdot \prod_j V_j^{-1}$. By the inequality
of the quadratic and geometric means this expression is maximal
if and only if $b_1=b_2=\dots =b_n$, which is achieved by the
choice $c_j:=V_j^{-1}$. Substituting $c_j:=V_j^{-1}$ in
(\ref{becs}) we arrive at the required inequality (\ref{opt}).
\end{proof}

We remark that $V_j\ge 1$ for every $j$, and $V_j= 1$ for
all $j$ if only $x_1,\dots,x_n$ forms an orthonormal system, hence in general we cannot prove the Conjecture by this
argument. It is easy to see, however, that our result is stronger
than the estimate \eqref{marcus} of Marcus.
\begin{corollary}
Assume that the given unit vectors $x_1,\dots , x_n$ are linearly independent.
Let $\l_1, \dots ,\l_n$ denote the eigenvalues of the Gram matrix $XX\adj $. Then
\begin{equation}\label{har}
\max_{\|y\|=1}\prod_{j=1}^{n}|(Xy)_j| \ge
\left(\frac{n}{\l_1^{-1}+ \cdots +\l_n^{-1}}\right)^{n/2}\cdot
n^{n/2}
\end{equation}
\end{corollary}

\begin{proof}
Observe that $\sum_k V_k^2=\sum_{jk}v_{jk}^2=\tr (XX\adj )^{-1}= \l_1^{-1}+ \cdots +\l_n^{-1}$.
Therefore
$$
\frac{1}{V_1\cdots V_n}\ge \left(\frac{n}{\sum_k
V_k^2}\right)^{n/2}=\left(\frac{n}{\l_1^{-1}+ \cdots
+\l_n^{-1}}\right)^{n/2}.
$$
\end{proof}

The next two estimates are in the same spirit but the proofs rely
on Bang's lemma instead of the averaging technique of Theorem
\ref{prop1}.

\begin{theorem}\label{prop2}
Let $\l_1, \dots ,\l_n$ denote the eigenvalues of the Gram matrix $XX\adj $ in increasing order. Then
\begin{equation}\label{nagy}
\max_{\|y\|=1}\prod_{j=1}^{n}|(Xy)_j|\ge
\left(\frac{1}{\l_n}\right)^{n/2}\cdot n^{n/2}
\end{equation}
\end{theorem}

\begin{proof}
Let $\e:=(\e_1,\dots \e_n)$ be any sequence of signs, and let $E$ denote the diagonal matrix with the numbers $(\e_j)_{j=1}^n$ in the diagonal. Bang's lemma (applied to the numbers $r_1=r_2=\dots =r_n=1$) implies that there exists a sequence of signs $\e$ such that $E XX\adj E (1,1, \dots ,1)^T\ge (1,1,\dots ,1)^T$. Take
$$y:=\frac{(XX\adj )\sq \e}{\| (XX\adj )\sq \e\|}.$$
Then
\begin{equation}
\prod_j |((XX\adj )\sq y)_j| \ge \left(\frac{1}{\|(XX\adj )\sq
\e\|}\right)^n \ge (n \l_n)^{-n/2}
\end{equation}
where the last estimate follows from $\|\e\|=n\sq$ and $\|(XX\adj )\sq \|=\l_n\sq$.

\end{proof}

We remark that the choice of $y$ above actually satisfies the
inequality $\prod_j |((XX\adj )\sq y)_j| \ge n^{-n/2}$ (proving
also the Conjecture) for $1\le n \le 5$, as shown in \cite{papp}.
It is not clear, however, whether this particular choice works
also for larger values of $n$. Note that $\l_n\ge 1$, and Theorem
\ref{prop2} yields the conjectured estimate only in case of
$\l_n= 1$, that is, only in the case of an orthonormal system.

\begin{theorem}\label{prop3}
Let $a_1, \dots ,a_n$ denote the diagonal entries of the matrix $(XX\adj )^{1/2}$. Then
\begin{equation}\label{diag}
\max_{\|y\|=1}\prod_{j=1}^{n}|(Xy)_j|
\ge {a_1 \cdots a_n}\cdot n^{-n/2} \, .
\end{equation}
\end{theorem}

\begin{proof}
Let $A$ denote the $n\times n$ diagonal matrix with
$(a_j)_{j=1}^n$ in the diagonal. The matrix $B:= A^{-1/2}(XX\adj
 )\sq A^{-1/2}$ is positive, self-adjoint and has 1's in the main
diagonal, i.e. it is a Gram matrix. Apply Bang's lemma to $B$
with numbers $a_1\sq , \dots , a_n\sq$. We conclude that there
exists a choice of signs $\e:=(\e_1, \dots , \e_n)$ (with the
corresponding diagonal matrix $E$) such that $EA\sq BA\sq E
(1,1,\dots ,1)^T\ge (a_1,a_2,\dots , a_n)^T$ (coordinatewise).
This means that the choice $y:= n^{-1/2}(\e_1,\dots ,\e_n)^T$
gives $\prod_j |((XX\adj )\sq y)_j|\ge {a_1 \cdots a_n}\cdot
n^{-\frac{n}{2}}$.
\end{proof}

Note that $a_j\ge \l_1\sq$ for every $j$, therefore the result
above is also an improvement on Marcus' estimate.

\medskip

Let us conclude with the following remarks.

The three different estimates presented in Theorems \ref{prop1},
\ref{prop2} and  \ref{prop3} do not seem to be comparable; at
least we were not able to show that any of them would imply
another. An advantage of the proofs applied is that in all three
theorems we were able to pinpoint our choice of the vector $y$.
It is clear, however, that none of these estimates settles the
Conjecture.

\end{document}